\documentclass[amstex,12pt,russian,amssymb]{article}

\usepackage{mathtext}
\usepackage[cp1251]{inputenc}
\usepackage[T2A]{fontenc}
\usepackage[russian]{babel}
\usepackage[dvips]{graphicx}
\usepackage{amsmath}
\usepackage{amssymb}
\usepackage{amsxtra}
\usepackage{latexsym}
\usepackage{ifthen}

\begin{document}

\newcounter{lemma}
\newcommand{\lemma}{\par \refstepcounter{lemma}%
{\bf Лема \arabic{lemma}.}}

\newcounter{corollary}
\newcommand{\corollary}{\par \refstepcounter{corollary}%
{\bf Наслідок \arabic{corollary}.}}

\newcounter{remark}
\newcommand{\remark}{\par \refstepcounter{remark}%
{\bf Зауваження \arabic{remark}.}}

\newcounter{theorem}
\newcommand{\theorem}{\par \refstepcounter{theorem}%
{\bf Теорема \arabic{theorem}.}}

\newcounter{proposition}
\newcommand{\proposition}{\par \refstepcounter{proposition}%
{\bf Твердження \arabic{proposition}.}}

\newcounter{example}
\newcommand{\example}{\par \refstepcounter{example}%
{\bf Приклад \arabic{example}.}}

\renewcommand{\refname}{\centerline{\bf Список літератури}}

\renewcommand{\figurename}{Мал.}

\newcommand{\proof}{{\it Доведення.\,\,}}

\noindent УДК 517.5

\medskip\medskip
{\bf О.П.~Довгопятий} (Житомирський державний університет імені
Івана Фран\-ка)

{\bf Є.О.~Севостьянов} (Житомирський державний університет імені
Івана Фран\-ка; Інститут прикладної математики і механіки НАН
України, м.~Слов'янськ)

\medskip\medskip
{\bf O.P.~Dovhopiatyi} (Zhytomyr Ivan Franko State University)

{\bf E.A.~Sevost'yanov} (Zhytomyr Ivan Franko State University;
Institute of Applied Ma\-the\-ma\-tics and Mechanics of NAS of
Ukraine, Slov'yans'k)

\medskip
{\bf Про відображення з гідродинамічним нормуванням у евклідовому
просторі}

{\bf On mappings with hydrodynamical normalization conditions in
Euclidean space}

\medskip\medskip
Вивчаються просторові відображення, які задовольняють деякий
просторовий аналог гідродинамічної умови зростання в околі
нескінченно віддаленої точки. Доведено, що гомеоморфізми вказаного
класу формують одностайно неперервні сім'ї за деяких умов на їх
характеристику квазіконформності. Розглянуто також питання щодо
замкненості цих класів відносно локально рівномірної збіжності.
Отримані відповідні результати для відображень з інтегральними
обмеженнями, а також для класів відповідних обернених відображень.

\medskip\medskip
We are studying spatial mappings that satisfy some space analog of a
hydrodynamical type of growth in the neighborhood of the infinity.
It is proved that homeomorphisms of the specified class form
equicontinuous families under some conditions on their
characteristic of quasiconformality. We have also considered the
problem of closeness of these classes with respect to locally
uniform convergence. We have obtained corresponding results for
mappings with integral constraints, as well as for classes of
corresponding inverse mappings.

\newpage
{\bf 1. Вступ.} Дана стаття присвячена відображенням з узагальненням
так званого гідродинамічного нормування $f(z)=z+o(1),$ $z\rightarrow
\infty.$ Їх вивчення є важливим, зокрема, з огляду на існування
відповідних гомеоморфних розв'язків рівняння Бельтрамі (див., напр.,
\cite[теорема~1.B.V]{A}, \cite[теореми~1.1--1.2]{GRSY$_1$},
\cite[лема~1]{GRSY$_2$}). В наших роботах вже досліджувалася
проблема компактності класів таких розв'язків
(див.~\cite{DS$_1$}--\cite{DS$_2$}). Дослідження відбувалися
виключно на площині, і основною ціллю даного рукопису є
розповсюдження аналогічних результатів у евклідовий $n$-вимірний
простір. Зауважимо, що <<просторовий>> випадок відрізняється від
<<плоского>> з огляду відсутності прямого аналогу теореми Кебе про
чверть (ця теорема істотно використовувалася при доведенні;
див.~\cite[теорема~1.3]{CG}). Іншою відмінністю є відсутність
просторових аналогів рівнянь Бельтрамі, тому основні результати
наводяться в термінах більш абстрактних класів відображень.

\medskip
Наведемо деякі означення. Нехай $D$ -- область в ${\Bbb R}^n,$
$n\geqslant 2,$ $M(\Gamma)$ -- конформний модуль сім'ї кривих
$\Gamma$ в ${\Bbb R}^n$ (див., напр., \cite[гл.~6]{Va}). Покладемо
$$S(x_0, r)=\{x\in {\Bbb R}^n: |x-x_0|=r\}\,,B(x_0, r)=\{x\in {\Bbb R}^n:
|x-x_0|<r\}\,,$$
$${\Bbb B}^n:=B(0, 1)\,,\quad {\Bbb S}^{n-1}:=S(0, 1)\,,\quad \Omega_n=m({\Bbb B}^n)\,,\quad
\omega_{n-1}=\mathcal{H}^{n-1}({\Bbb S}^{n-1})\,,$$
де $\mathcal{H}^{n-1}$ позначає $(n-1)$-вимірну міру Хаусдорфа в
${\Bbb R}^n.$ Нехай, крім того,
$$A=A(x_0, r_1, r_2)=\{x\in {\Bbb R}^n: r_1<|x-x_0|<r_2\}\,.$$
Для заданих множин $E,$ $F\subset\overline{{\Bbb R}^n}$ і області
$D\subset {\Bbb R}^n$ позначимо через $\Gamma(E,F,D)$ сім'ю всіх
кривих $\gamma:[a,b]\rightarrow \overline{{\Bbb R}^n}$ таких, що
$\gamma(a)\in E,\gamma(b)\in\,F$ і $\gamma(t)\in D$ при $t \in [a,
b].$ Відображення $f:D\rightarrow \overline{{\Bbb R}^n}$ називається
{\it кільцевим $Q$-ві\-доб\-ра\-жен\-ням у точці
$x_0\,\in\,\overline{D},$} якщо співвідношення
\begin{equation}\label{eq3*!gl0}
M(f(\Gamma(S(x_0, r_1),\,S(x_0, r_2),\,D)))\leqslant
\int\limits_{A\cap D} Q(x)\cdot \eta^n(|x-x_0|)\, dm(x)
\end{equation}
виконано для будь-якого кільця $A=A(x_0, r_1,r_2),$ $0<r_1<r_2<
r_0:=\sup\limits_{x\in D}|x-x_0|,$ і кожної вимірної за Лебегом
функції $\eta:(r_1,r_2)\rightarrow [0,\infty ]$ такої, що
\begin{equation}\label{eq*3gl0}
\int\limits_{r_1}^{r_2}\eta(r)\,dr \geqslant 1\,.
\end{equation}

\medskip
Нехай $D$ -- область в ${\Bbb R}^n.$ Будемо говорити, що функція
${\varphi}:D\rightarrow{\Bbb R},$ що є локально інтегровною в
деякому околі точки $x_0\in D,$ має {\it скінченне середнє
коливання} в точці $x_0$ (пишемо: $\varphi\in FMO(x_0)$), якщо
\begin{equation}\label{eq17:}
{\limsup\limits_{\varepsilon\rightarrow
0}}\frac{1}{\Omega_n\varepsilon^n}\int\limits_{B(
x_0,\,\varepsilon)}
|{\varphi}(x)-\overline{{\varphi}}_{\varepsilon}|\ dm(x)\, <\,
\infty\,,
\end{equation}
де $\Omega_n$ -- об'єм одиничної кулі в ${\Bbb R}^n,$
$\overline{{\varphi}}_{\varepsilon}=\frac{1}{\Omega_n\varepsilon^n}\int\limits_{B(
x_0,\,\varepsilon)} {\varphi}(x)\ dm(x)$ (див., напр.,
\cite[розд.~2]{RSY$_2$}).

\medskip Нехай $Q:{\Bbb R}^n\rightarrow [0, \infty]$ -- вимірна за Лебегом
функція і $K$ -- компакт у ${\Bbb R}^n.$ Позначимо через
$\frak{F}_{Q}(K)$ клас усіх гомеоморфізмів $f:{\Bbb R}^n\rightarrow
{\Bbb R}^n,$ які задовольняють умову~(\ref{eq3*!gl0}) в кожній точці
$x_0\in {\Bbb R}^n,$ $f(x)\ne 0$ при всіх $x\not\in K,$ $f\in
ACL({\Bbb R}^n\setminus K),$ причому для будь-якого $r_0>0$
знайдеться число $M_0=M_0(r_0)>0$ таке, що
\begin{equation}\label{eq8}
a_{\widetilde{f}_{r_0}}(0):=\exp\left(\frac{1}{n\Omega_n}
\int\limits_{{\Bbb B}^n}\log |J(x, \widetilde{f}_{r_0})|\,dm(x)
\right)\geqslant M_0\,,
\end{equation}
де $\widetilde{f}_{r_0}(x):=\frac{1}{r_0}\cdot(\psi\circ f\circ
\psi)(r_0x)$ і $\psi(x)=\frac{x}{|x|^2}.$

\medskip Одним з основних результатів рукопису є наступне твердження.
\medskip
\begin{theorem}\label{th1}
{\sl\, Нехай функція $Q$ задовольняє умову $Q(x)\leqslant Q_0=const$
при майже всіх $x\in {\Bbb R}^n\setminus K,$ крім того, припустимо,
що $Q$ задовольняє принаймні одну з умов:

\medskip
1) або $Q\in FMO({\Bbb R}^n),$

\medskip
2) або для кожного $x_0\in {\Bbb R}^n$ існує
$\delta_0=\delta(x_0)>0$ таке, що
\begin{equation}\label{eq2D}
\int\limits_0^{\delta_0}\frac{dt}{tq^{\frac{1}{n-1}}_{x_0}(t)}=\infty\,,
\end{equation}
де $q_{x_0}(t)=\frac{1}{\omega_{n-1}t^{n-1}}\int\limits_{S(x_0,
t)}Q(x)\,d\mathcal{H}^{n-1}.$
Тоді сім'я відображень $\frak{F}_{Q}(K)$ є одностайно неперервною в
${\Bbb R}^n.$}
\end{theorem}

\medskip
В теоремі~\ref{th1} одностайну неперервність слід розуміти як між
метричними просторами $(X, d)$ і $\left(X^{\,\prime},
d^{\,\prime}\right),$ де $X={\Bbb R}^n,$
$X^{\,\prime}=\overline{{\Bbb R}^n},$ $d$ -- евклідова метрика і
$d^{\,\prime}=h$ -- хордальна (сферична метрика),
$$h(x,\infty)=\frac{1}{\sqrt{1+{|x|}^2}}\,,$$
\begin{equation}\label{eq3C}
\ \ h(x,y)=\frac{|x-y|}{\sqrt{1+{|x|}^2} \sqrt{1+{|y|}^2}}\,, \ \
x\ne \infty\ne y
\end{equation}
(див., напр., \cite[означення~12.1]{Va}).

\medskip
\begin{remark}\label{rem1}
Умова~(\ref{eq8}) в означенні класу $\frak{F}_{Q}(K)$ є істотною,
оскільки доволі легко побудувати приклад сім'ї $\frak{F}$
гомеоморфізмів $f:{\Bbb R}^n\rightarrow {\Bbb R}^n$ з $Q\equiv 1$ у
${\Bbb R}^n$ таких, щоб виконувалася умова~(\ref{eq3*!gl0}) та умови
1)--2) з формулювання теореми~\ref{th1}; проте, в той самий час, ця
сім'я не була одностайно неперервною в ${\Bbb R}^n.$ В якості такої
сім'ї можна взяти, наприклад, клас відображень $f_m(x)=mx,$
$m=1,2,\ldots .$
\end{remark}

\medskip
{\bf 2. Основна лема.} Наступний просторовий аналог теореми Кебе про
чверть належить Астала та Герінгу, див. \cite[теорема~1.8]{AG}.

\medskip
\begin{proposition}\label{pr1}
{\sl\,Припустимо, що $D$ і $D^{\,\prime}$ -- області в ${\Bbb R}^n,$
$n\geqslant 2.$ Якщо $f:D\rightarrow D^{\,\prime}$ --
$K$-квазіконформне відображення, то для всіх $x\in D$ виконуються
нерівності
\begin{equation}\label{eq9}
\frac{1}{c}\cdot\frac{d(f(x), \partial D^{\,\prime})}{d(x, \partial
D)} \leqslant a_f(x)\leqslant c\cdot \frac{d(f(x), \partial
D^{\,\prime})}{d(x, \partial D)}\,,
\end{equation}
де $c$ -- стала, залежна тільки від $K$ і $n,$ та
$$a_{f}(x):=\exp\left(\frac{1}{n\Omega_nd^n(x,
\partial D)} \int\limits_{B(x, d(x,
\partial D))}\log |J(x, \widetilde{f})|\,dm(x) \right)\,.$$
 }
\end{proposition}

\medskip
\begin{remark}\label{rem2}
Як зазначено в~\cite[теорема~1.6]{AG}, для конформних відображень
площини ($n=2$) виконуються рівності: $a_f(x)=|f^{\,\prime}(x)|$ і
$c=4.$ Якщо $D={\Bbb D}$ -- одиничний круг, $x=0=f(0)$ і
$f^{\,\prime}(0)=1,$ то з правої частини нерівності в~(\ref{eq9}) ми
отримаємо, що $1\leqslant 4\cdot d(0, \partial D^{\,\prime}),$ або
$d(0, \partial D^{\,\prime})\geqslant \frac{1}{4}.$ Звідси випливає,
що $B(0, 1/4)\subset f({\Bbb D}),$ що є змістом класичної теореми
Кебе.

З міркувань, наведених у доведенні теореми~1.1 в~\cite{DS$_2$},
випливає, що з умови $f(z)=z+o(1),$ $z\rightarrow\infty,$ $z\in
{\Bbb C},$ виконаної для відображення $f,$ конформного в деякому
околі нескінченності, випливає, що $\widetilde{f}(0)=0$ і
$\widetilde{f}^{\,\prime}(0)=1$ (де $\widetilde{f}$ визначено нижче
формули~(\ref{eq8})). Отже, умова~(\ref{eq8}) в цьому випадку може
бути замінена умовою $f(z)=z+o(1),$ $z\rightarrow\infty.$
\end{remark}

\medskip
Справедливе наступне твердження.

\medskip
\begin{lemma}\label{lem1}
{\sl\, Нехай функція $Q$ задовольняє умову $Q(x)\leqslant Q_0=const$
при майже всіх $x\in {\Bbb R}^n\setminus K.$ Крім того, припустимо,
що існує число $\varepsilon_0>0$ та невід'ємна вимірна за Лебегом
функція $\psi:(0,\infty)\rightarrow [0,\infty]$ така, що у кожній
точці $x_0\in K$ виконується умова
\begin{equation} \label{eq4!}
\int\limits_{A(x_0, \varepsilon,
\varepsilon_0)}Q(x)\cdot\psi^n(|x-x_0|) \
dm(x)\,=\,o\left(I^n(\varepsilon, \varepsilon_0)\right)\,,
\end{equation}
де
\begin{equation} \label{eq5A}
0< I(\varepsilon,
\varepsilon_0):=\int\limits_{\varepsilon}^{\varepsilon_0}\psi(t)dt <
\infty\qquad\forall\,\, \varepsilon\in (0, \varepsilon_0)\,,
\end{equation}
причому $I(\varepsilon, \varepsilon_0)\rightarrow \infty$ при
$\varepsilon\rightarrow 0.$ Тоді сім'я відображень $\frak{F}_{Q}(K)$
є одностайно неперервною в ${\Bbb R}^n.$}
\end{lemma}

\medskip
\begin{proof} Зафіксуємо $f\in\frak{F}_{Q}(K)$ і $r_0>0$ -- число, таке що $K\subset B(0, 1/r_0).$
Покладемо
\begin{equation}\label{eq1A}
\widetilde{f}_{r_0}(x)=\frac{1}{r_0}\cdot\frac{f\left(\frac{x}{r_0|x|^2}\right)}
{\left|f\left(\frac{x}{r_0|x|^2}\right)\right|^2}\,\,,\qquad x\in
{\Bbb B}^n\setminus\{0\}\,.
\end{equation}
Зауважимо, що для кожного відображення $f\in \frak{F}_{Q}(K)$
виконується умова
\begin{equation}\label{eq1B}
\Vert f^{\,\prime}(x)\Vert^n \leqslant C_n\cdot|J(x,f)|\ Q^{n-1}(x)
\end{equation}
майже скрізь, де $\Vert
f^{\,\prime}(x)\Vert=\sup\limits_{|h|=1}|f^{\,\prime}(x)h|,$ $J(x,
f)=\det f^{\,\prime}(x)$ і $C_n>0$ -- деяка стала, залежна тільки
від розмірності простору $n$ (див., напр.,
\cite[наслідок~3.4]{SalSev}). Отже, з нерівності~(\ref{eq1B})
випливає, що при майже всіх $x\in {\Bbb R}^n\setminus K$
$$\Vert f^{\,\prime}(x)\Vert^n \leqslant C_n\cdot|J(x,f)|\ Q_0^{n-1}\,,$$
отже,
\begin{equation}\label{eq1C}
K_O(x, f)\leqslant C_n\cdot Q_0^{n-1}<\infty
\end{equation}
майже скрізь, де зовнішня дилатація відображення $f$ у точці $x$
обчислюється за правилом
\begin{equation}\label{eq18}
K_O(x,f)\quad =\quad \left\{
\begin{array}{rr}
\frac{\Vert f^{\,\prime}(x)\Vert^n}{|J(x,f)|}, & J(x,f)\ne 0,\\
1,  &  f^{\,\prime}(x)=0, \\
\infty, & {\rm в\,\,інших\,\,випадках}
\end{array}
\right.\,.
\end{equation}
Зауважимо, що $f$ є диференційовним майже скрізь (див., напр.,
\cite[теорема~3.2]{SalSev}). Тоді оскільки за припущенням $f\in
ACL({\Bbb R}^n\setminus K),$ з огляду на нерівність~(\ref{eq1C})
відображення $f$ є $K_0$-квазіконформним в ${\Bbb R}^n\setminus K,$
де $K_0$ -- деяке число, залежне тільки від розмірності простору $n$
(див., напр., \cite[теорема~34.6]{Va}).

\medskip
В такому випадку, зауважимо, що таким є і відображення
$\widetilde{f}_{r_0}$ визначене по $f$ у~(\ref{eq1A}). Дійсно,
$\widetilde{f}_{r_0}=\psi_{r_0}\circ f\circ \psi_{r_0},$ де
$\psi_{r_0}(x)=\frac{1}{r_0}\cdot\frac{x}{|x|^2},$ причому, оскільки
$\psi_{r_0}$ є конформним відображенням і
$\widetilde{f}_{r_0}^{\,\prime}(x)=\psi_{r_0}^{\,\prime}(f(\psi_{r_0}(x)))\cdot
f^{\,\prime}(\psi_{r_0}(x))\cdot \psi_{r_0}^{\,\prime}(x),$ то з
огляду на співвідношення~(\ref{eq1C}) та обчислення дилатацій від
суперпозиціями з конформними відображеннями (див. \cite[розд.~4,
гл.~I]{Re}), будемо мати, що
$$K_O(x, \widetilde{f}_{r_0})=K_O(f(\psi_{r_0}(x)), \psi_{r_0})\cdot K_O(\psi_{r_0}(x),
f)\cdot K_0(x, \psi_{r_0})=$$
\begin{equation}\label{eq2B}
=1\cdot K_O(\psi_{r_0}(x), f)\cdot 1=K_O(\psi_{r_0}(x), f)\leqslant
C_n\cdot Q_0^{n-1}
\end{equation}
для всіх $x\in {\Bbb B}^n.$ Очевидно, відображення $f$ є
гомеоморфізмом в ${\Bbb B^n}\setminus\{0\},$ і що воно є
диференційовним майже скрізь. Більше того, оскільки конформне
відображення $\psi_{r_0}(x)=\frac{1}{r_0}\cdot\frac{x}{|x|^2}$ є
локальною квазіізометрією, то як внутрішня, так і зовнішня
суперпозиція відображення $f$ з ним не виводить за межі класу $ACL$
(див., напр., \cite[теорема~1, розд.~1.1.7 гл.~1]{Ma}). Тоді знову
$\widetilde{f}_{r_0}$ є $K_0$-квазіконформним в ${\Bbb
B}^n\setminus\{0\},$ де $K_0$ -- деяке число, залежне тільки від
розмірності простору $n$ (див., напр., \cite[теорема~34.6]{Va}).
Зауважимо, що $\widetilde{f}_{r_0}$ має квазіконформне продовження в
точку $x_0=0$ (див., напр., \cite[теорема~17.3]{Va}). В такому
випадку, само відображення $f$ має квазіконформне продовження в
точку $\infty.$

\medskip
Зауважимо, що $\widetilde{f}_{r_0}(0)=0.$ Дійсно, $f(\overline{{\Bbb
R}^n})$ є одночасно відкритою і замкненою підмножиною
$\overline{{\Bbb R}^n}$ як, з одного боку, відкритий образ відкритої
множини $\overline{{\Bbb R}^n},$ а з іншого -- неперервний образ
компакту $\overline{{\Bbb R}^n}.$ Отже, $f(\overline{{\Bbb
R}^n})=\overline{{\Bbb R}^n}.$ Оскільки $f(x)\ne \infty$ у ${\Bbb
R}^n,$ то $f(\infty)=\infty.$ Звідси випливає, що
$\widetilde{f}_{r_0}(0)=0,$ що і було потрібно.

\medskip
Позначимо
\begin{equation}\label{eq1AA}
\frak{R}_{Q}(K):=\left\{ \widetilde{f}_{r_0}:{\Bbb B}^n\rightarrow
{\Bbb R}^n:
\widetilde{f}_{r_0}(x)=\frac{1}{r_0}\cdot\frac{f\left(\frac{x}{r_0|x|^2}\right)}
{\left|f\left(\frac{x}{r_0|x|^2}\right)\right|^2}, f\in
\frak{F}_{Q}(K)\right\}\,.
\end{equation}

\medskip
Нехай $\widetilde{f}_{r_0}\in\frak{R}_{Q}(K).$ З огляду на
співвідношення~(\ref{eq8}) $a_{\widetilde{f}_{r_0}}(0)\geqslant
M_0=M_0(r_0).$ Тоді за твердженням~\ref{pr1}
\begin{equation}\label{eq4}
B(0, c/M_0)\subset \widetilde{f}_{r_0}({\Bbb B}^n)\,.
\end{equation}
Зі співвідношення~(\ref{eq4}) випливає, що
\begin{equation}\label{eq5}
F(\overline{{\Bbb R}^n}\setminus \overline{B(0, 1/r_0)})\supset B(0,
cr_0/M_0)\,,\qquad F(x)=f(x)/|f(x)|^2\,.
\end{equation}
З урахуванням~(\ref{eq5}) покажемо, що
\begin{equation}\label{eq6}
f(\overline{{\Bbb R}^n}\setminus \overline{B(0, 1/r_0)})\supset
\overline{{\Bbb R}^n}\setminus \overline{B(0, M_0/cr_0)}\,.
\end{equation}
Дійсно, нехай $y\in \overline{{\Bbb R}^n}\setminus \overline{B(0,
M_0/cr_0)},$ тоді $\frac{y}{|y|^2}\in B(0, cr_0/M_0).$ Зі
співвідношення~(\ref{eq5}) будемо мати, що
$\frac{y}{|y|^2}=f(x)/|f(x)|^2,$ $x\in \overline{{\Bbb
R}^n}\setminus \overline{B(0, 1/r_0)}.$ Тоді $y=f(x),$ $x\in
\overline{{\Bbb R}^n}\setminus \overline{B(0, 1/r_0)},$ що і
доводить~(\ref{eq6}).

\medskip
Оскільки $f$ -- гомеоморфізм у ${\Bbb R}^n,$ зі
співвідношення~(\ref{eq6}) випливає, що
\begin{equation}\label{eq10}
f(B(0, 1/r_0))\subset B(0, M_0/cr_0)\,.
\end{equation}
Покладемо $\Delta:=h(\overline{{\Bbb R}^n}\setminus B(0,
M_0/cr_0)),$ де $h(\overline{{\Bbb R}^n}\setminus B(0, M_0/cr_0))$
-- хордальний діаметр множини $\overline{{\Bbb R}^n}\setminus B(0,
M_0/cr_0).$ В таком випадку, сім'я відображень $\frak{F}_Q(K)$ є
одностайно неперервною в~$B(0, 1/r_0)$ за~\cite[лема~7.6]{MRSY}.
Отже, $\frak{F}_Q(K)$ є одностайно неперервною в ${\Bbb R}^n,$ бо
число $r_0>0$ було обрано довільним додатнім.~$\Box$
\end{proof}

\medskip
{\it Твердження теореми~\ref{th1}} безпосередньо випливає з
леми~\ref{lem1}, бо виконання умов 1) і 2) в цій теоремі є частковим
випадком співвідношень~(\ref{eq4!})--(\ref{eq5A}) зі спеціально
обраними функціями $\psi.$ Деталі див., напр., в \cite[наслідок~6.3,
теорема~6.4]{MRSY}.~$\Box$

\medskip
{\bf 3. Замкненість одного підкласу $\frak{F}_{Q}(K)$.} В повній
мірі питання щодо замкненості класу $\frak{F}_{Q}(K)$ залишається
відкритим, оскільки невідомо, чи можливий граничний перехід у
нерівностях типу~(\ref{eq8}). Тим не менш, справедливим є наведене
нижче твердження.

\medskip Нехай $Q:{\Bbb R}^n\rightarrow [0, \infty]$ -- вимірна за Лебегом
функція і $K$ -- компакт у ${\Bbb R}^n.$ Позначимо через
$\frak{M}_{Q}(K)$ клас усіх гомеоморфізмів $f:{\Bbb R}^n\rightarrow
{\Bbb R}^n,$ які задовольняють умову~(\ref{eq3*!gl0}) в кожній точці
$x_0\in {\Bbb R}^n,$ $f(x)\ne 0$ при всіх $x\not\in K,$ $f\in
ACL({\Bbb R}^n\setminus K),$ причому для будь-якого $r_0>0$
знайдуться числа $M_0=M_0(r_0)>0$ і $N_0=N_0(r_0)>0$ такі, що
\begin{equation}\label{eq10A}
N_0\geqslant
a_{\widetilde{f}_{r_0}}(0):=\exp\left(\frac{1}{n\Omega_n}
\int\limits_{{\Bbb B}^n}\log |J(x, \widetilde{f}_{r_0})|\,dm(x)
\right)\geqslant M_0\,,
\end{equation}
де $\widetilde{f}_{r_0}(x):=\frac{1}{r_0}\cdot(\psi\circ f\circ
\psi)(r_0x)$ і $\psi(x)=\frac{x}{|x|^2}.$ Справедлива наступна

\medskip
\begin{lemma}\label{lem2}
{\sl\, Нехай функція $Q$ задовольняє умову $Q(x)\leqslant Q_0=const$
при майже всіх $x\in {\Bbb R}^n\setminus K,$ крім того, припустимо,
що виконуються умови~(\ref{eq4!})--(\ref{eq5A}).

Нехай, крім того, $f_m\in \frak{M}_{Q}(K),$ $m=1,2,\ldots ,$ --
послідовність відображень класу $\frak{M}_{Q}(K),$ яка збігається до
деякого відображення $f:{\Bbb R}^n\rightarrow \overline{{\Bbb R}^n}$
локально рівномірно у ${\Bbb R}^n$ при $m\rightarrow\infty$ відносно
хордальної метрики $h.$ Тоді $f(x)\ne \infty$ при всіх $x\in {\Bbb
R}^n,$ крім того, $f$ є гомеоморфізмом, який задовольняє
умову~(\ref{eq3*!gl0}) в кожній точці $x_0\in {\Bbb R}^n,$ $f(x)\ne
0$ при всіх $x\not\in K,$ $f\in ACL({\Bbb R}^n\setminus K),$ причому
для будь-якого $r_0>0$ знайдуться числа $M_*=M_*(r_0, f)>0$ і
$N_*=N_*(r_0, f)>0$ такі, що виконується умова
\begin{equation}\label{eq9A}
N_*\geqslant
a_{\widetilde{f}_{r_0}}(0):=\exp\left(\frac{1}{n\Omega_n}
\int\limits_{{\Bbb B}^n}\log |J(x, \widetilde{f}_{r_0})|\,dm(x)
\right)\geqslant M_*\,,
\end{equation}
де $\widetilde{f}_{r_0}(x):=\frac{1}{r_0}\cdot(\psi\circ f\circ
\psi)(r_0x)$ і $\psi(x)=\frac{x}{|x|^2}.$  }
\end{lemma}

\medskip
\begin{remark}
Зауважимо, що в лемі~\ref{lem2} числа $M_*$ і $N_*$ можуть, взагалі
кажучи, залежати від відображення~$f.$
\end{remark}

\medskip
{\it Доведення леми~\ref{lem2}}. Нехай $f_m,$ $m=1,2,\ldots ,$ --
відображення з умови леми. За~\cite[лема~4.2]{RSS} відображення $f$
є або гомеоморфізмом $f:{\Bbb R}^n\rightarrow {\Bbb R}^n,$ або
сталою $c_0\in\overline{{\Bbb R}^n}.$

\medskip
Покажемо, що друга ситуація неможлива. Насамперед, міркуючи
аналогічно доведенню леми~\ref{lem1} та з огляду на
оцінку~(\ref{eq10A}), ми отримаємо, що
\begin{equation}\label{eq11}
f_m(B(0, 1/r_0))\subset B(0, M_0/r_0)\,.
\end{equation}
Співвідношення~(\ref{eq11}) виключає випадок $c_0=\infty.$

\medskip
Нехай тепер $c_0\ne\infty.$ Розглядаючи сім'ю відображень
$\widetilde{f_m}_{r_0}$ по аналогії з~(\ref{eq1A}),
\begin{equation}\label{eq14}
\widetilde{f_m}_{r_0}(x)=\frac{1}{r_0}\cdot\frac{f_m\left(\frac{x}{r_0|x|^2}\right)}
{\left|f_m\left(\frac{x}{r_0|x|^2}\right)\right|^2}\,\,,\qquad x\in
{\Bbb B}^n\,.
\end{equation}
з огляду на ліву нерівність у~(\ref{eq9}) та на оцінку~(\ref{eq10A})
будемо мати, що
$$d(0, \partial \widetilde{f_m}_{r_0}({\Bbb
B}^n))\leqslant c a_{\widetilde{f}_{r_0}}(0)\leqslant cN_*\,,\quad
m=1,2,\ldots \,.$$
Тоді
\begin{equation}\label{eq11A}
B(0, 2cN_*)\cap \left(\overline{{\Bbb
R}^n}\setminus\widetilde{f_m}_{r_0}({\Bbb B}^n)\right)\ne
\varnothing\,,\quad m=1,2,\ldots \,.
\end{equation}
З~(\ref{eq11A}) випливає, що існує $x_m\in B(0, 2cN_*)$ такий, що
$x_m\in B(0, 2cN_*)\cap \left(\overline{{\Bbb
R}^n}\setminus\widetilde{f_m}_{r_0}({\Bbb B}^n)\right).$ Крім того,
$\widetilde{f_m}_{r_0}(x)\ne \infty$ при всіх $x\in {\Bbb B}^n$ за
означенням класу $\frak{M}_{Q}(K)$ та числа $r_0.$ Тоді
\begin{equation}\label{eq12}
h\left(\overline{{\Bbb R}^n}\setminus\widetilde{f_m}_{r_0}({\Bbb
B}^n)\right)\geqslant h(x_m,
\infty)=\frac{1}{\sqrt{1+|x_m|^2}}\geqslant
\frac{1}{\sqrt{1+|2cN_*|^2}}:=\delta_0\,.
\end{equation}
Оскільки всі відображення $\widetilde{f_m}$ квазіконформні з
загальною сталою квазіконформності, $m=1,2,\ldots $ (див. доведення
леми~\ref{lem1}), то при додатковій умові~(\ref{eq12}) ця сім'я є
одностайно неперервною (див.~\cite[теорема~19.2]{Va}). Зокрема, для
числа $\varepsilon:=\frac{1}{r_0(c_0+1)}$ знайдеться
$\delta=\delta(\varepsilon)>0$ таке, що
\begin{equation}\label{eq13}
|\widetilde{f_m}_{r_0}(x)|<\frac{1}{r_0(c_0+1)}\qquad\forall\quad
|x|\leqslant \delta, \,\,m=1,2,\ldots\,.
\end{equation}
Тоді з огляду на~(\ref{eq13}) відповідно до~(\ref{eq14}) ми
отримаємо, що $\frac{1}{r_0}\biggl|\cdot\frac{f_m\left(x\right)}
{\left|f_m\left(x\right)\right|^2}\biggr|<\frac{1}{r_0(c_0+1)}$ при
$|x|\geqslant \frac{1}{r_0\delta},$ або
\begin{equation}\label{eq15}
|f_m(x)|>c_0+1\,,\qquad |x|\geqslant \frac{1}{r_0\delta}\,,\quad
m=1,2,\ldots .
\end{equation}
Якщо ж тепер $f_m(x)\rightarrow c_0=const$ при $m\rightarrow\infty$
локально рівномірно в ${\Bbb R}^n$ при $m\rightarrow\infty,$ то для
числа $A=1$ знайдеться номер $M^*=M^*(r_0, \delta)$ такий, що
\begin{equation}\label{eq16}
|f_m(x)-c_0|<1\,,\quad |x|\leqslant \frac{2}{r_0\delta}\,,\quad
m\geqslant M^*\,.
\end{equation}
З~(\ref{eq16}) та з нерівності трикутника випливає, що
\begin{equation}\label{eq17}
|f_m(x)|\leqslant c_0+1\,,\quad |x|\leqslant
\frac{2}{r_0\delta}\,,\quad m\geqslant M^*\,.
\end{equation}
Нерівності~(\ref{eq15}) та~(\ref{eq17}) суперечать одна одній, отже,
відображення $f$ є або гомеоморфізмом $f:{\Bbb R}^n\rightarrow {\Bbb
R}^n,$ що і треба було довести. Далі, з огляду на
умови~(\ref{eq4!})--(\ref{eq5A}) відображення $f$ задовольняє також
визначальні співвідношення~(\ref{eq3*!gl0})--(\ref{eq*3gl0})
(див.~\cite[теорема~5.1]{RSS}).

\medskip
Далі, $f$ є квазіконформним відображенням на кожній зв'язній
компоненті відкритої множини ${\Bbb R}^n\setminus K$ як локально
рівномірна границя квазіконформних відображень
(див.~\cite[наслідок~37.3]{Va}). Отже, $f\in ACL({\Bbb R}^n\setminus
K)$ (див.~\cite[наслідок~31.4]{Va}).

\medskip
Доведемо, що $f(x)\ne 0$ при $x\in {\Bbb R}^n\setminus K.$ Дійсно,
оскільки $f_m({\Bbb R}^n)={\Bbb R}^n$ і $f_m(x)\ne 0$ при $x\in
{\Bbb R}^n\setminus K,$ то існує $x_m\in K$ таке, що $f_m(x_m)=0.$
Оскільки $K$ -- компакт, то можна виділити підпослідовність
$x_{m_k}\in K,$ $k=1,2,\ldots $ таку, що $x_{m_k}\rightarrow x_0,$
$x_0\in K.$ Тоді за нерівністю трикутника та з огляду на локально
рівномірну збіжність $f_{m_k}$ до $f$ будемо мати, що
$$|f_{m_k}(x_{m_k})-f(x_0)|\leqslant |f_{m_k}(x_{m_k})-f_{m_k}(x_0)|+
|f_{m_k}(x_0)-f(x_0)|\rightarrow 0$$
при $k\rightarrow\infty.$ Оскільки $f_{m_k}(x_{m_k})=0,$ останнє
співвідношення є можливим тільки за умови $f(x_0)=0.$ Оскільки $f$ є
гомеоморфізмом у ${\Bbb R}^n,$ такою точкою, в якій обертається в
нуль відображення $f,$ може бути лише $x_0.$ Тому $f(x)\ne 0$ при
$x\in {\Bbb R}^n\setminus K,$ що і потрібно було довести.

\medskip
Нарешті, залишилося встановити співвідношення~(\ref{eq9A}). Оскільки
відображення $\widetilde{f}_{r_0}$ є квазіконформним (див. деталі
доведення леми~\ref{lem1}), то це співвідношення виконується з
огляду на твердження~\ref{pr1}. Лема~\ref{lem2} повністю
доведена.~$\Box$

\medskip
Справедлива наступна

\medskip
\begin{theorem}\label{th2}
{\sl\, Нехай функція $Q$ задовольняє умову $Q(x)\leqslant Q_0=const$
при майже всіх $x\in {\Bbb R}^n\setminus K,$ крім того, припустимо,
що $Q\in FMO({\Bbb R}^n),$ або для кожного $x_0\in {\Bbb R}^n$ існує
$\delta_0=\delta(x_0)>0$ таке, що виконується умова~(\ref{eq2D}).
Нехай, крім того, $f_m\in \frak{M}_{Q}(K),$ $m=1,2,\ldots ,$ --
послідовність відображень класу $\frak{M}_{Q}(K),$ яка збігається до
деякого відображення $f:{\Bbb R}^n\rightarrow \overline{{\Bbb R}^n}$
локально рівномірно у ${\Bbb R}^n$ при $m\rightarrow\infty$ відносно
хордальної метрики $h.$ Тоді $f(x)\ne \infty$ при всіх $x\in {\Bbb
R}^n,$ крім того, $f$ є гомеоморфізмом, який задовольняє
умову~(\ref{eq3*!gl0}) в кожній точці $x_0\in {\Bbb R}^n,$ $f(x)\ne
0$ при всіх $x\not\in K,$ $f\in ACL({\Bbb R}^n\setminus K),$ причому
для будь-якого $r_0>0$ знайдуться числа $M_*=M_*(r_0, f)>0$ і
$N_*=N_*(r_0, f)>0$ такі, що виконується умова~(\ref{eq9A}), де
$\widetilde{f}_{r_0}(x):=\frac{1}{r_0}\cdot(\psi\circ f\circ
\psi)(r_0x)$ і $\psi(x)=\frac{x}{|x|^2}.$}
\end{theorem}

\medskip
{\it Доведення теореми~\ref{th2}} безпосередньо випливає з
леми~\ref{lem2}, бо вказані у першій частині формулювання умови на
функцію $Q$ є окремим випадком
співвідношень~(\ref{eq4!})--(\ref{eq5A}) зі спеціально обраними
функціями $\psi$ (див., напр., \cite[наслідок~6.3,
теорема~6.4]{MRSY}).~$\Box$

\medskip
{\bf 4. Про класи відображень з інтегральними обмеженнями.} Нехай
$\Phi\colon\overline{{\Bbb R}^{+}}\rightarrow\overline{{\Bbb
R}^{+}}$ -- неспадна функція, $Q_0, L_0>0$ -- фіксовані числа і $K$
-- компакт у ${\Bbb R}^n.$ Позначимо через $\frak{F}^{Q_0,
L_0}_{\Phi}(K)$ клас усіх гомеоморфізмів $f:{\Bbb R}^n\rightarrow
{\Bbb R}^n,$ для яких:

\medskip
\noindent 1) знайдеться вимірна за Лебегом функція $Q=Q_f:{\Bbb
C}\rightarrow[0, \infty]$ така, що $Q(x)\leqslant Q_0=const$ майже
скрізь у ${\Bbb R}^n\setminus K;$

\medskip
\noindent 2) виконана умова~(\ref{eq3*!gl0}) в кожній точці $x_0\in
{\Bbb R}^n,$ причому
\begin{equation}\label{e3.3.1}
\int\limits_{{\Bbb
R}^n}\Phi(Q(x))\cdot\frac{dm(x)}{(1+|x|^2)^n}\leqslant L_0\,;
\end{equation}

\medskip
\noindent 3) $f(x)\ne 0$ при всіх $x\not\in K,$ $f\in ACL({\Bbb
R}^n\setminus K);$

\medskip
\noindent 4) для будь-якого $r_0>0$ знайдеться число
$M_0=M_0(r_0)>0$ таке, що
\begin{equation}\label{eq8A}
a_{\widetilde{f}_{r_0}}(0):=\exp\left(\frac{1}{n\Omega_n}
\int\limits_{{\Bbb B}^n}\log |J(x, \widetilde{f}_{r_0})|\,dm(x)
\right)\geqslant M_0\,,
\end{equation}
де $\widetilde{f}_{r_0}(x):=\frac{1}{r_0}\cdot(\psi\circ f\circ
\psi)(r_0x)$ і $\psi(x)=\frac{x}{|x|^2}.$

\medskip
Виконується наступний результат.

\medskip
\begin{theorem}\label{th3}
{\sl\, Нехай $\Phi:\overline{{\Bbb R^{+}}}\rightarrow
\overline{{\Bbb R^{+}}}$ -- непрерывна зростаюча опукла функція, яка
при деякому $\delta>\Phi(0)$ задовольняє умову
\begin{equation}\label{eq2A}
\int\limits_{\delta}^{\infty}\frac{d\tau}{\tau\left(\Phi^{\,-1}(\tau)\right)^{\frac{1}{n-1}}}=\infty\,.
\end{equation}
Тоді сім'я відображень $\frak{F}^{Q_0, L_0}_{\Phi}(K)$ є одностайно
неперервною в ${\Bbb R}^n.$}
\end{theorem}

\medskip
{\it Доведення теореми~\ref{th3}} дуже схоже на доведення
леми~\ref{lem1}, тому обмежимося лише схематичним доведенням.
Зафіксуємо $f\in\frak{F}^{Q_0, L_0}_{\Phi}(K)$ і довільний компакт
$C\subset {\Bbb R}^n.$ Нехай $r_0>0$ -- число, таке що $C\cup
K\subset B(0, 1/r_0).$ Визначимо $\widetilde{f}_{r_0}$ за
співвідношенням~(\ref{eq1A}). Аналогічно доведенню леми~\ref{lem1}
можна показати, що $f$ є $K_0$-квазіконформним в ${\Bbb
R}^n\setminus K,$ де $K_0$ -- деяке число, залежне тільки від
розмірності простору $n.$ Так само, як і в лемі~\ref{lem1} можна
показати, що і відображення $\widetilde{f}_{r_0}$ є
$K_0$-квазіконформним в ${\Bbb B}^n,$ причому
$\widetilde{f}_{r_0}(0)=0.$

\medskip
Позначимо
\begin{equation}\label{eq1AAA}
\frak{R}^{Q_0, L_0}_{\Phi}(K):=\left\{ \widetilde{f}_{r_0}:{\Bbb
B}^n\rightarrow {\Bbb R}^n:
\widetilde{f}_{r_0}(x)=\frac{1}{r_0}\cdot\frac{f\left(\frac{x}{r_0|x|^2}\right)}
{\left|f\left(\frac{x}{r_0|x|^2}\right)\right|^2}, f\in
\frak{F}^{Q_0, L_0}_{\Phi}(K)\right\}\,.
\end{equation}

\medskip
Нехай $\widetilde{f}_{r_0}\in \frak{R}^{Q_0, L_0}_{\Phi}(K).$ З
огляду на співвідношення~(\ref{eq8}) будемо мати, що
\begin{equation}\label{eq6A}
f(\overline{{\Bbb R}^n}\setminus \overline{B(0, 1/r_0)})\supset
\overline{{\Bbb R}^n}\setminus \overline{B(0, M_0/cr_0)}\,.
\end{equation}
Оскільки $f$ -- гомеоморфізм у ${\Bbb R}^n,$ зі
співвідношення~(\ref{eq6A}) випливає, що
\begin{equation}\label{eq10AA}
f(B(0, 1/r_0))\subset B(0, M_0/cr_0)\,.
\end{equation}
Покладемо $\Delta:=h(\overline{{\Bbb R}^n}\setminus B(0,
M_0/cr_0)),$ де $h(\overline{{\Bbb R}^n}\setminus B(0, M_0/cr_0))$
-- хордальний діаметр множини $\overline{{\Bbb R}^n}\setminus B(0,
M_0/cr_0).$ В таком випадку, сім'я відображень $\frak{F}^{Q_0,
L_0}_{\Phi}(K)$ є одностайно неперервною в~$B(0, 1/r_0)$
за~\cite[теорема~4.1]{RS}.~$\Box$

\medskip
Нехай $\Phi\colon\overline{{\Bbb R}^{+}}\rightarrow\overline{{\Bbb
R}^{+}}$ -- неспадна функція, $Q_0, L_0>0$ -- фіксовані числа і $K$
-- компакт у ${\Bbb R}^n.$ Позначимо через $\frak{M}^{Q_0,
L_0}_{\Phi}(K)$ клас усіх гомеоморфізмів $f:{\Bbb R}^n\rightarrow
{\Bbb R}^n,$ для яких:

\medskip
\noindent 1) знайдеться вимірна за Лебегом функція $Q=Q_f:{\Bbb
C}\rightarrow[0, \infty]$ така, що $Q(x)\leqslant Q_0=const$ майже
скрізь у ${\Bbb R}^n\setminus K;$

\medskip
\noindent 2) виконана умова~(\ref{eq3*!gl0}) в кожній точці $x_0\in
{\Bbb R}^n,$ причому
\begin{equation}\label{e3.3.1A}
\int\limits_{{\Bbb
R}^n}\Phi(Q(x))\cdot\frac{dm(x)}{(1+|x|^2)^n}\leqslant L_0\,;
\end{equation}

\medskip
\noindent 3) $f(x)\ne 0$ при всіх $x\not\in K,$ $f\in ACL({\Bbb
R}^n\setminus K);$

\medskip
\noindent 4) для будь-якого $r_0>0$ знайдуться числа
$N_0=N_0(r_0)>0$ і $M_0=M_0(r_0)>0$ таке, що
\begin{equation}\label{eq8B}
N_0\geqslant
a_{\widetilde{f}_{r_0}}(0):=\exp\left(\frac{1}{n\Omega_n}
\int\limits_{{\Bbb B}^n}\log |J(x, \widetilde{f}_{r_0})|\,dm(x)
\right)\geqslant M_0\,,
\end{equation}
де $\widetilde{f}_{r_0}(x):=\frac{1}{r_0}\cdot(\psi\circ f\circ
\psi)(r_0x)$ і $\psi(x)=\frac{x}{|x|^2}.$

\medskip
Виконується наступний результат.

\medskip
\begin{theorem}\label{th4}
{\sl\, Нехай $\Phi:\overline{{\Bbb R^{+}}}\rightarrow
\overline{{\Bbb R^{+}}}$ -- непрерывна зростаюча опукла функція, яка
при деякому $\delta>\Phi(0)$ задовольняє умову
\begin{equation}\label{eq2AA}
\int\limits_{\delta}^{\infty}\frac{d\tau}{\tau\left(\Phi^{\,-1}(\tau)\right)^{\frac{1}{n-1}}}=\infty\,.
\end{equation}
Нехай, крім того, $f_m\in \frak{M}^{Q_0, L_0}_{\Phi}(K),$
$m=1,2,\ldots ,$ -- послідовність відображень класу $\frak{M}^{Q_0,
L_0}_{\Phi}(K),$ яка збігається до деякого відображення $f:{\Bbb
R}^n\rightarrow \overline{{\Bbb R}^n}$ локально рівномірно у ${\Bbb
R}^n$ при $m\rightarrow\infty$ відносно хордальної метрики $h.$ Тоді
$f(x)\ne \infty$ при всіх $x\in {\Bbb R}^n,$ крім того, $f$ є
гомеоморфізмом, який задовольняє умову~$f(x)\ne 0$ при всіх
$x\not\in K,$ $f\in ACL({\Bbb R}^n\setminus K),$ причому для
будь-якого $r_0>0$ знайдуться числа $M_*=M_*(r_0, f)>0$ і
$N_*=N_*(r_0, f)>0$ такі, що виконується умова
\begin{equation}\label{eq9AA}
N_*\geqslant
a_{\widetilde{f}_{r_0}}(0):=\exp\left(\frac{1}{n\Omega_n}
\int\limits_{{\Bbb B}^n}\log |J(x, \widetilde{f}_{r_0})|\,dm(x)
\right)\geqslant M_*\,,
\end{equation}
де $\widetilde{f}_{r_0}(x):=\frac{1}{r_0}\cdot(\psi\circ f\circ
\psi)(r_0x)$ і $\psi(x)=\frac{x}{|x|^2}.$

}
\end{theorem}

\medskip
{\it Доведення теореми~\ref{th4}} дуже схоже на доведення
леми~\ref{lem2}. Отже, обмежимося лише схемою доведення.

\medskip
Нехай $f_m,$ $m=1,2,\ldots ,$ -- відображення з умови теореми.
За~\cite[лема~2.1]{DS$_2$} відображення $f$ є або гомеоморфізмом
$f:{\Bbb R}^n\rightarrow {\Bbb R}^n,$ або сталою
$c_0\in\overline{{\Bbb R}^n}.$

\medskip
Покажемо, що друга ситуація неможлива. Міркуючи аналогічно доведенню
леми~\ref{lem1} та з огляду на оцінку~(\ref{eq8B}), ми отримаємо, що
\begin{equation}\label{eq11B}
f_m(B(0, 1/r_0))\subset B(0, M_0/r_0)\,.
\end{equation}
Співвідношення~(\ref{eq11B}) виключає випадок $c_0=\infty.$

\medskip
Нехай тепер $c_0\ne\infty.$ Розглядаючи сім'ю відображень
$\widetilde{f_m}_{r_0}$ по аналогії з~(\ref{eq1A}),
\begin{equation}\label{eq14A}
\widetilde{f_m}_{r_0}(x)=\frac{1}{r_0}\cdot\frac{f_m\left(\frac{x}{r_0|x|^2}\right)}
{\left|f_m\left(\frac{x}{r_0|x|^2}\right)\right|^2}\,\,,\qquad x\in
{\Bbb B}^n\,,
\end{equation}
і міркуючи аналогічно доведенню леми~\ref{lem2}, ми отримаємо, що
\begin{equation}\label{eq12A}
h\left(\overline{{\Bbb R}^n}\setminus\widetilde{f_m}_{r_0}({\Bbb
B}^n)\right)\geqslant\delta_0
\end{equation}
для деякого числа $\delta_0>0.$ Оскільки всі відображення
$\widetilde{f_m}$ квазіконформні з загальною сталою
квазіконформності, $m=1,2,\ldots $ (див. доведення леми~\ref{lem1}),
то при додатковій умові~(\ref{eq12A}) ця сім'я є одностайно
неперервною (див.~\cite[теорема~19.2]{Va}). Зокрема, для числа
$\varepsilon:=\frac{1}{r_0(c_0+1)}$ знайдеться
$\delta=\delta(\varepsilon)>0$ таке, що
\begin{equation}\label{eq13A}
|\widetilde{f_m}_{r_0}(x)|<\frac{1}{r_0(c_0+1)}\qquad\forall\quad
|x|\leqslant \delta, \,\,m=1,2,\ldots\,.
\end{equation}
Тоді з огляду на~(\ref{eq13A}) відповідно до~(\ref{eq14A}) ми
отримаємо, що $\frac{1}{r_0}\biggl|\cdot\frac{f_m\left(x\right)}
{\left|f_m\left(x\right)\right|^2}\biggr|<\frac{1}{r_0(c_0+1)}$ при
$|x|\geqslant \frac{1}{r_0\delta},$ або
\begin{equation}\label{eq15A}
|f_m(x)|>c_0+1\,,\qquad |x|\geqslant \frac{1}{r_0\delta}\,,\quad
m=1,2,\ldots .
\end{equation}
Якщо ж тепер $f_m(x)\rightarrow c_0=const$ то (міркуючи аналогічно
доведенню леми~\ref{lem2}) ми отримаємо, що
\begin{equation}\label{eq17A}
|f_m(x)|\leqslant c_0+1\,,\quad |x|\leqslant
\frac{2}{r_0\delta}\,,\quad m\geqslant M^*
\end{equation}
для деякого числа $M^*>0.$ Нерівності~(\ref{eq15A}) та~(\ref{eq17A})
суперечать одна одній, отже, відображення $f$ є або гомеоморфізмом
$f:{\Bbb R}^n\rightarrow {\Bbb R}^n,$ що і треба було довести.

\medskip
Те, що $f\in ACL({\Bbb R}^n\setminus K)$ і $f(x)\ne 0$ при $x\in
{\Bbb R}^n\setminus K$ доводиться так само, як у лемі~\ref{lem2}, як
і співвідношення~(\ref{eq9AA}).~$\Box$

\medskip
{\bf 5. Відображення з оберненою нерівністю Полецького.} Одностайна
неперервність таких відображень при різних умовах досліджувалася в
багатьох наших роботах (див., напр., \cite{DS$_2$} і \cite{SSD}).
Нижче наводиться ще один тип умов в контексті цього питання.

\medskip
Нагадаємо означення. Якщо $f:D\rightarrow {\Bbb R}^n$ -- задане
відображення, $y_0\in f(D)$ і $0<r_1<r_2<d_0=\sup\limits_{y\in
f(D)}|y-y_0|,$ то через $\Gamma_f(y_0, r_1, r_2)$ ми позначимо сім'ю
всіх кривих $\gamma$ в області $D$ таких, що $f(\gamma)\in
\Gamma(S(y_0, r_1), S(y_0, r_2), A(y_0,r_1,r_2)).$ Нехай $Q:{\Bbb
R}^n\rightarrow [0, \infty]$ -- вимірна за Лебегом функція.  Будемо
говорити, що {\it $f$ задовольняє обернену нерівність Полецького} в
точці $y_0\in f(D),$ якщо співвідношення
\begin{equation}\label{eq2*A}
M(\Gamma_f(y_0, r_1, r_2))\leqslant \int\limits_{A(y_0,r_1,r_2)\cap
f(D)} Q(y)\cdot \eta^n (|y-y_0|)\, dm(y)
\end{equation}
виконується для довільної вимірної за Лебегом функції $\eta:
(r_1,r_2)\rightarrow [0,\infty ]$ такій, що
\begin{equation}\label{eqA2}
\int\limits_{r_1}^{r_2}\eta(r)\, dr\geqslant 1\,.
\end{equation}

\medskip Нехай $Q:{\Bbb R}^n\rightarrow [0, \infty]$ -- вимірна за Лебегом
функція і $K$ -- компакт у ${\Bbb R}^n.$ Позначимо через
$\frak{B}_{Q}(K)$ клас усіх гомеоморфізмів $f:{\Bbb R}^n\rightarrow
{\Bbb R}^n,$ які є гомеоморфізмами у ${\Bbb R}^n\setminus K$ і
задовольняють умови~(\ref{eq2*A})--(\ref{eqA2}) в кожній точці
$x_0\in {\Bbb R}^n,$ $f(x)\ne 0$ при всіх $x\not\in K,$
$g=f^{\,-1}\in ACL(f({\Bbb R}^n\setminus K)),$ причому для
будь-якого $r_0>0$ знайдеться число $M_0=M_0(r_0)>0$ таке, що
\begin{equation}\label{eq8C}
a_{\widetilde{f}_{r_0}}(0):=\exp\left(\frac{1}{n\Omega_n}
\int\limits_{{\Bbb B}^n}\log |J(x, \widetilde{f}_{r_0})|\,dm(x)
\right)\geqslant M_0\,,
\end{equation}
де $\widetilde{f}_{r_0}(x):=\frac{1}{r_0}\cdot(\psi\circ f\circ
\psi)(r_0x)$ і $\psi(x)=\frac{x}{|x|^2}.$ Справедлива наступна

\medskip
\begin{theorem}\label{th1}
{\sl\, Нехай функція $Q$ задовольняє умову $Q(x)\leqslant Q_0=const$
при майже всіх $x\in {\Bbb R}^n\setminus K,$ крім того, припустимо,
що для кожної точки $y_0\in {\Bbb R}^n$ і кожних $0<r_1<r_2<\infty$
знайдеться множина $E\subset[r_1, r_2]$ додатної міри Лебега така,
що при кожному $r\in E$ функція $Q$ є інтегровною по відношенню до
$\mathcal{H}^{n-1}$-міри на сфері $S(y_0, r).$ Тоді сім'я
відображень $\frak{F}_{Q}(K)$ є одностайно неперервною в ${\Bbb
R}^n.$}
\end{theorem}

\medskip
\begin{proof}
Зафіксуємо $f\in\frak{B}_{Q}(K)$ і $r_0>0$ -- число, таке що $
K\subset B(0, 1/r_0).$ Нехай, крім того, $\widetilde{f}_{r_0}$
визначається формулою~(\ref{eq1A}). З умови~(\ref{eq2*A}) випливає,
що співвідношення~(\ref{eq3*!gl0}) виконується в $f({\Bbb
R}^n\setminus K)$ для відображення $g:=f^{\,-1}$ при $Q\mapsto Q_0.$

\medskip
Зауважимо, що для кожного відображення $g=f^{\,-1},$ $f\in
\frak{B}_{Q}(K),$ виконується умова
\begin{equation}\label{eq1BB}
\Vert g^{\,\prime}(y)\Vert^n \leqslant C_n\cdot|J(y,f)|\ Q^{n-1}(y)
\end{equation}
майже скрізь, де $\Vert
g^{\,\prime}(x)\Vert=\sup\limits_{|h|=1}|g^{\,\prime}(x)h|,$ $g(x,
f)=\det f^{\,\prime}(x)$ і $C_n>0$ -- деяка стала, залежна тільки
від розмірності простору $n$ (див., напр.,
\cite[наслідок~3.4]{SalSev}). Отже, з нерівності~(\ref{eq1BB})
випливає, що при майже всіх $y\in f({\Bbb R}^n\setminus K)$
$$\Vert g^{\,\prime}(x)\Vert^n \leqslant C_n\cdot|J(x,g)|\ Q_0^{n-1}\,,$$
отже,
\begin{equation}\label{eq1C}
K_O(y, g)\leqslant C_n\cdot Q_0^{n-1}<\infty
\end{equation}
майже скрізь, де зовнішня дилатація $K_O(y, g)$ відображення $g$ у
точці $y$ обчислюється за правилом~(\ref{eq18}). Зауважимо, що $g$ є
диференційовним майже скрізь в $f({\Bbb R}^n\setminus K)$ (див.,
напр., \cite[теорема~3.2]{SalSev}). Тоді оскільки за припущенням
$g\in ACL(f({\Bbb R}^n\setminus K)),$ з огляду на
нерівність~(\ref{eq1C}) відображення $g$ є $K_0$-квазіконформним в
$f({\Bbb R}^n\setminus K),$ де $K_0$ -- деяке число, залежне тільки
від розмірності простору $n$ (див., напр., \cite[теорема~34.6]{Va}).
Тоді також і відображення $f$ є $K_0$-квазіконформним в ${\Bbb
R}^n\setminus K$ (див. \cite[означення~13.1]{Va}).

\medskip
В такому випадку, зауважимо, що таким є і відображення
$\widetilde{f}_{r_0}$ визначене по $f$ у~(\ref{eq1A}), що
встановлюється так само, як і при доведенні леми~\ref{lem1}.
Зауважимо, що $\widetilde{f}_{r_0}$ має квазіконформне продовження в
точку $x_0=0$ (див., напр., \cite[теорема~17.3]{Va}). В такому
випадку, само відображення $f$ має квазіконформне продовження в
точку $\infty.$

\medskip
Зауважимо, що $\widetilde{f}_{r_0}(0)=0$ (це також встановлюється
так, як при доведенні леми~\ref{lem1}). Кріс того, міркуючи
аналогічно доведенню леми~\ref{lem1}, ми отримаємо співвідношення
\begin{equation}\label{eq10B}
f(B(0, 1/r_0))\subset B(0, M_0/cr_0)
\end{equation}
для кожного $f\in\frak{B}_{Q}(K).$ В такому випадку, сім'я
відображень є одностайно неперервною в $B(0, 1/r_0)$
(див.~\cite[теорема~1.1]{SevSkv$_1$}). Отже, $\frak{F}_Q(K)$ є
одностайно неперервною в ${\Bbb R}^n,$ бо число $r_0>0$ було обрано
довільним додатнім.~$\Box$
\end{proof}

\medskip
\begin{corollary}\label{cor1}
{\sl\,Твердження теореми~\ref{th4} залишається справедливим, якщо в
ньому замість відповідних умов на $Q$ вимагати, щоб $Q(x)\leqslant
Q_0=const$ при майже всіх $x\in {\Bbb R}^n\setminus K$ і крім того,
що $Q\in L^1_{\rm loc}({\Bbb R}^n).$ В цьому випадку, для будь-якої
компактної множини $C$ у ${\Bbb R}^n$ і області $D\subset {\Bbb
R}^n$ такої, що $C\subset D$ знайдеться область $D^{\,\prime}\subset
{\Bbb R}^n$ і функція $Q^{\,\prime},$ яка дорівнює $Q$ у
$D^{\,\prime}$ і обертається в нуль у її доповненні така що
нерівність
\begin{equation}\label{eq2E}
|f(x)-f(y)|\leqslant\frac{C_0}{\log^{1/n}\left(1+\frac{\widetilde{r_0}}{2|x-y|}\right)}
\end{equation}
виконується для будь-яких $x, y\in C$ і всіх $f\in
f\in\frak{B}_{Q}(K),$ де $C_0=C_0(n, C, \Vert Q^{\,\prime}\Vert_1,
D, D^{\,\prime})>0$ -- деяка стала, залежна тільки від $n,$ $C$ і
$\Vert Q^{\,\prime}\Vert_1,$ $\Vert Q^{\,\prime}\Vert_1$ позначає
$L^1$-норму функції $Q^{\,\prime}$ в ${\Bbb R}^n,$ крім того,
$\widetilde{r_0}=d(C,
\partial D).$}
\end{corollary}

\medskip
\begin{proof} За теоремою Фубіні (див.~\cite[теорема~8.1.III]{Sa})
для будь-якого $\varepsilon_0>0$ і $y_0\in {\Bbb R}^n$
$$\int\limits_{B(y_0, \varepsilon_0)}Q(y)\,dm(y)=\int\limits_{0}^{\varepsilon_0}
\int\limits_{S(y_0,
t)}Q(y)\,d\mathcal{H}^{n-1}(y)\,dt=\omega_{n-1}\int\limits_{0}^{\varepsilon_0}t^{n-1}
q_{y_0}(t)\,dt<\infty\,.$$
Отже, на майже всіх сферах $S(y_0, r)$ функція $Q$ інтегровна, тому
функція $Q$ задовольняє всі умови теореми~\ref{th4}. Повторюючи
доведення цієї теореми, приходимо до співвідношення~(\ref{eq10B}).
Якщо тепер маємо довільний компакт $C\subset {\Bbb R}^n,$ то завжди
можна підібрати $r_0>0$ таким, що $C\subset B(0, 1/r_0).$ Далі
покладемо $D=B(0, 1/r_0)$ і $D^{\,\prime}=B(0, M_0/cr_0).$ В такому
випадку, співвідношення~(\ref{eq2E}) випливає
з~\cite[теорема~4.1]{SevSkv$_1$}.~$\Box$
\end{proof}

%=================Список литературы====================
%\end{fulltext}

КОНТАКТНА ІНФОРМАЦІЯ

\medskip
\noindent{{\bf Олександр Петрович Довгопятий} \\
Житомирський державний університет ім.\ І.~Франко\\
вул. Велика Бердичівська, 40 \\
м.~Житомир, Україна, 10 008 \\
e-mail: alexdov1111111@gmail.com}

\medskip
\noindent{{\bf Євген Олександрович Севостьянов} \\
{\bf 1.} Житомирський державний університет ім.\ І.~Франко\\
вул. Велика Бердичівська, 40 \\
м.~Житомир, Україна, 10 008 \\
{\bf 2.} Інститут прикладної математики і механіки
НАН України, \\
вул.~Добровольського, 1 \\
м.~Слов'янськ, Україна, 84 100\\
e-mail: esevostyanov2009@gmail.com}


\begin{thebibliography}{99}

{\small


\bibitem{A} {\it Альфорс Л.} Лекции по квазиконформным
отображениям. -- Москва: Мир, 1969.

\bibitem{GRSY$_1$} {\it Gutlyanskii~V., Ryazanov~V., Sevost'yanov~E., Yakubov~E.}
On the degenerate Beltrami equation and hydrodynamic normalization
// Journal of Mathematical Sciences. -- 2022. -- \textbf{262,}
no.~2. -- P.~165--183.

\bibitem{GRSY$_2$} {\it Gutlyanskii~V., Ryazanov~V., Sevost'yanov~E., Yakubov~E.}
BMO and Dirichlet problem for degenerate Beltrami equation //
Journal of Mathematical Sciences. -- 2022. -- \textbf{268,} no.~2.
-- P.~155--177.

\bibitem{DS$_1$} {\it Dovhopiatyi~O.P., Sevost'yanov~E.A.}
On the compactness of classes of the solutions of the Dirichlet
problem // Journal of Mathematical Sciences. -- 2021. --
\textbf{259}, no.~1. -- P.~23--36.

\bibitem{DS$_2$} {\it Dovhopiatyi~O.P., Sevost'yanov~E.A.} On compact classes of
Beltrami solutions and Dirichlet problem // Complex Variables and
Elliptic Equations. -- 2022. --
https://www.tandfonline.com/doi/abs/10.1080/17476933.2022.2040020 .

\bibitem{CG} {\it Carleson, L., Gamelin~T.W.} Complex dynamics. -
Universitext: Tracts in Mathematics, New York etc.: Springer-Verlag,
1993.

\bibitem{Va} {\it V\"{a}is\"{a}l\"{a}~J.} Lectures on $n$-Dimensional Quasiconformal
Mappings. --  Lecture Notes in Math.~\textbf{229,} Berlin etc.:
Springer--Verlag, 1971.

\bibitem{RSY$_2$} {\it Ryazanov~V., Srebro~U. and Yakubov~E.}
Finite mean oscillation and the Beltrami equation // Israel Math. J.
-- 2006. -- \textbf{153.} -- P.~247--266.

\bibitem{AG} {\it Astala K., Gehring F.} Quasiconformal analogues of theorems of Koebe and
Hardy–Littlewood // Mich. Math. J.  -- 1985. -- \textbf{32,} no.~1.
-- P.~99–-107.

\bibitem{SalSev} {\it Salimov~R.R. and Sevost'yanov~E.A.} ACL and differentiability of
the open discrete ring mappings // Complex Variables and Elliptic
Equations. -- 2010. -- \textbf{55,} no.~1-3. -- P.~49--59.

\bibitem{Re} {\it Reshetnyak Yu.G.} Space Mappings with Bounded Distortion.
Transl. of Math. Monographs~\textbf{73}, AMS, 1989.

\bibitem{Ma} {\it Maz'ya V.} Sobolev classes. -- New York:
Springer, Berlin, 1985.

\bibitem{MRSY} {\it Martio O.,  Ryazanov V.,  Srebro U.,
 Yakubov E.} Moduli in modern mapping theory. -- New York: Springer Science + Business Media, LLC,
 2009.

\bibitem{RSS} {\it Ryazanov~V., Salimov~R. and Sevost'yanov~E.}
On Convergence Analysis of Space Homeomorphisms // Siberian Advances
in Mathematics. -- 2013. -- \textbf{23}, no.~4. -- P.~263--293.

\bibitem{RS} {\it Ryazanov~V., Sevost'yanov~E.}
Equicontinuity of mappings quasiconformal in the mean // Ann. Acad.
Sci. Fenn. -- 2011. -- \textbf{36.} -- P. 231--244.

\bibitem{SSD} {\it Sevost’yanov E.A., Skvortsov~S.A. and
Dovhopiatyi~O.P.} On nonhomeomorphic mappings with the inverse
Poletsky inequality // Journal of Mathematical Sciences. -- 2021. --
\textbf{252,} no.~4. -- P.~541--557.

\bibitem{SevSkv$_1$} {\it Sevost’yanov E.A. and Skvortsov~S.A.}
Logarithmic H\"{o}lder continuous mappings and Beltrami equation //
Analysis and Mathematical Physics. -- 2021. -- \textbf{11}, No. 3.
-- Article number 138.

\bibitem{Sa} {\it Saks S.} Theory of the Integral. -- New York: Dover, 1964.
}


\end{thebibliography}
\end{document}